\documentclass[notitlepage]{amsart}

\usepackage{amsfonts}
\usepackage[latin1]{inputenc}
\usepackage[all]{xy} 
\usepackage{latexsym,amssymb,amscd}
\usepackage{amsmath}
\newtheorem{pro}{Proposition}[section]
\newtheorem{teo}[pro]{Theorem}
\newtheorem{defi}[pro]{Definition}
\newtheorem{lem}[pro]{Lemma}
\newtheorem{cor}[pro]{Corollary}
\newtheorem{rk}[pro]{Remark}

\newtheorem{ex}[pro]{Example}

\newcommand{\la}{\left\langle}
\newcommand{\ra}{\right\rangle}

\newcommand{\F}{\mathfrak{F}}

\newcommand{\X}{\mathcal{X}}

\newcommand{\modu}{\mathrm{mod}}

\newcommand{\proj}{\mathrm{proj}}
\newcommand{\Tor}{\mathrm{Tor}}
\newcommand{\Ext}{\mathrm{Ext}}
\newcommand{\Hom}{\mathrm{Hom}}
\newcommand{\End}{\mathrm{End}}
\newcommand{\Ker}{\mathrm{Ker}}
\newcommand{\Coker}{\mathrm{Coker}}
\newcommand{\Ima}{\mathrm{Im}}

\newcommand{\id}{\mathrm{id}}
\newcommand{\Tr}{\mathrm{Tr}}
\newcommand{\rad}{\mathrm{rad}}

\newcommand{\add}{\mathrm{add}}

\newenvironment{dem}{\noindent\bf Proof. \rm }{$\ \Box$}
\usepackage{latexsym,amssymb,amscd}
\usepackage{amsmath}

\begin{document}
\title[split algebras and stratifying systems]{Split-by-nilpotent extensions algebras and stratifying systems}
\author{M. Lanzilotta, O. Mendoza, C. S\'aenz.}
\thanks{2010 {\it{Mathematics Subject Classification}}. Primary 16G10. Secondary 18G99.\\
The authors thank the financial support received from
Project PAPIIT-UNAM IN100810.}
\date{}
\begin{abstract} Let $\Gamma$ and $\Lambda$ be artin algebras such that $\Gamma$ is a split-by-nilpotent extension of $\Lambda$ by a two sided ideal $I$ of $\Gamma.$ Consider the so-called change of rings functors $G:={}_\Gamma\Gamma_\Lambda\otimes_\Lambda -$ and $F:={}_\Lambda
\Lambda_\Gamma\otimes_\Gamma -.$  In this paper, we find the necessary and sufficient conditions under which a stratifying system $(\Theta,\leq)$ in $\modu\Lambda$ can be 
lifted to a stratifying system $(G\Theta,\leq)$ in $\modu\,(\Gamma).$ Furthermore, by using the 
functors $F$ and $G,$ we study the relationship between their filtered categories of modules; and 
some connections with their corresponding standardly stratified algebras are stated (see Theorem \ref{GbarraInducido}, Theorem \ref{Gequiv} and Theorem \ref{ApliGss-alg}). Finally, a sufficient condition is given for stratifying systems in $\modu\,(\Gamma)$ in such a way that they can be restricted, through the functor $F,$ to stratifying systems in $\modu\,(\Lambda).$  
\end{abstract}
\maketitle
\section{Introduction.} 

Stratifying systems where introduced in \cite{ES,MMS1,MMS2,MSX,Webb} and developed in \cite{LMMS, MMSZ,MS,Li,MPV} with some applications, for example, in \cite{DT,E,ESch,H,HLM,HP,MMS3,M}.
\

Split-by-nilpotent extension algebras have been recently studied in various settings. For 
example, in almost split sequences \cite{AZ2}, tilting modules \cite{AM} and quasi-tilted, laura, 
shod and weakly-shod algebras \cite{AZ}. We study these extension algebras from the point of 
view of the theory of stratifying systems.

The paper is organized as follows. After a brief section of preliminaries, we devote Section 3 to 
the study of the functors 
$F={}_\Lambda\Lambda_\Gamma\otimes_\Gamma:\modu\,(\Gamma)\to \modu\,(\Lambda)$ and 
$G={}_\Gamma\Gamma_\Lambda\otimes_\Lambda -:\modu\,(\Lambda)\to \modu\,(\Gamma),$ where 
$\Gamma$ is a split-by-nilpotent extension of $\Lambda$ by a two sided ideal $I$ of $\Gamma.$ \\
In Section 4, we show that, for any $M\in\modu\,(\Lambda),$ the algebra $\End_\Gamma(GM)$ is always an split-by-nilpotent extension of $\End_\Lambda(M)$ by $\Hom_\Lambda(M,I\otimes M)$ (see 
Theorem \ref{SplitAlgG}).\\
The section 5 is the main section in the paper. We give necessary and sufficient conditions such 
that the image under $G,$ of a stratifying system in $\modu\,(\Lambda),$ is a stratifying system 
in $\modu\,(\Gamma).$ Here, the main results are \ref{GbarraInducido}, \ref{Gequiv},  
\ref{ApliGss-alg} and \ref{ApliGss-algSS}. Finally, in Section 6, we give a sufficient condition 
(see \ref{TeoRestF}) for stratifying systems in $\modu\,(\Gamma)$ in such a way that they can be 
restricted, through the functor $F,$ to stratifying systems in $\modu\,(\Lambda).$ 
 
\section{Preliminaries.}

Throughout this paper, the term algebra means {\it artin algebra} over a commutative artin 
ring $R.$ For an algebra $\Lambda,$ the category of finitely generated left $\Lambda$-modules is 
denoted by $\modu\,(\Lambda).$ We denote by $\proj\,(\Lambda)$ the full subcategory of 
$\modu\,(\Lambda)$ whose objects are the projective $\Lambda$-modules. Unless otherwise specified, all the modules are finitely generated. 
Furthermore, for any positive integer $t,$ we set $[1,t]:=\{1,2,\cdots,t\}.$
\

\begin{defi}\cite{ES,MMS1} Let $\Lambda$ be an algebra. A stratifying system $(\Theta,\leq),$ of size $t$ in $\modu\,(\Lambda),$ consists
of a family of indecomposable $\Lambda$-modules $\Theta=\{\Theta(i)\}_{i=1}^t$ and a linear order $\leq$ on the
set $[1,t],$ satisfying the following two conditions.
\begin{itemize}
 \item[(a)] $\Hom_\Lambda(\Theta(i),\Theta(j))=0$ if $i>j.$
 \item[(b)] $\Ext_\Lambda^1(\Theta(i),\Theta(j))=0$ if $i\geq j.$
\end{itemize}
\end{defi}

For a set $\Theta$ of $\Lambda$-modules, let $\F(\Theta)$ be the subcategory of $\modu\,(\Lambda)$ consisting of the
$\Lambda$-modules $M$ having a $\Theta$-filtration, that is,
a sequence of submodules $0=M_0\subseteq M_1\subseteq\cdots\subseteq M_s=M$ such that each factor 
$M_{i+1}/M_i$ is isomorphic to a module in $\Theta$ for all $i$.

\begin{defi}  \cite{MMS2} \label{definicion de epps}
Let $\Lambda$ be an algebra. An Ext-projective
 stratifying system $(\Theta,\underline{Q},\leq),$ of size $t$ in $\modu\,(\Lambda),$ consists
of two families of non-zero $\Lambda$-modules
$\Theta=\{\Theta(i)\}_{i=1}^t$ and $\underline{Q}=\{ Q(i)\}_{i=1}^t,$
with $Q(i)$ indecomposable for all $i$, and a linear order $\leq$ on the
set $[1,t],$ satisfying the following three conditions.
\begin{itemize}
 \item[(a)] $\Hom_\Lambda(\Theta(i),\Theta(j))=0$ if $i>j.$
 \item[(b)] For each $i\in[1,t],$ there is an exact sequence
$$\varepsilon_i: 0\longrightarrow K(i)\longrightarrow Q(i)\overset{\beta_i}{\longrightarrow}\Theta(i)\longrightarrow0,$$
with $K(i)\in \F(\{\Theta(j):j>i\}).$
 \item[(c)] $Q$ is $\Theta$-projective. That is $\Ext_\Lambda^1(Q,\mathbf{\Theta})=0,$ where $Q:=\oplus_{i=1}^t\,Q(i)$ and $\mathbf{\Theta}:=\oplus_{i=1}^t\,\Theta(i).$
\end{itemize}
\end{defi}
 
Recall that (see \cite[Corollary 2.13]{MMS2}) an Ext-projective stratifying system $(\Theta,\underline{Q},\leq)$ gives the stratifying system 
$(\Theta,\leq).$ Furthermore, for a given a stratifying system $(\Theta,\leq),$ we know by \cite[Corollary 2.15]{MMS2} that there is a unique, up 
to isomorphism, Ext-projective stratifying system $(\Theta,\underline{Q},\leq).$ So, it is said that $(\Theta,\underline{Q},\leq)$ is the 
Ext-projective stratifying system associated to the stratifying system $(\Theta,\leq).$ We also have the dual notion of the Ext-injective 
stratifying system $(\Theta,\underline{Y},\leq)$ associated to the stratifying system $(\Theta,\leq)$ (see \cite{ES,MMS1,MMS2}).
\

The stratifying systems are related with the so-called standardly stratified algebras and so we introduce this notion. Let $\Lambda$ be an algebra.  For $M,N\in\modu\,(\Lambda)$, the {\bf trace} $\Tr_M\,(N)$ of $M$ in $N$, 
is the $\Lambda$-submodule of $N$ generated by the images of all
morphisms from $M$ to $N$. 
\

We next recall the definition (see \cite{R,DR,ADL,Dlab}) of the class of standard $\Lambda$-modules. Let $n$ be the rank of the Grothendieck group 
$K_0\,(\Lambda)$. We fix a linear order $\leq$ on the set $[1,n]$ and a representative set ${}_\Lambda P=\{{}_{\Lambda}P(i)\;:\;i\in[1,n]\},$ 
containing one module of each iso-class of indecomposable projective $\Lambda$-modules. Observe, 
that the set ${}_\Lambda P$ determines the representative set ${}_\Lambda S=\{{}_{\Lambda}S(i)\}_{i=1}^n$ of simple $\Lambda$-modules, where  ${}_{\Lambda}S(i):={}_{\Lambda}P(i)/\rad\,({}_{\Lambda}P(i))$ for each $i.$
\

The set of {\bf standard $\Lambda$-modules} is 
${}_\Lambda\Delta=\{{}_\Lambda\Delta(i):i\in[1,n]\},$ where
${}_\Lambda\Delta(i)={}_{\Lambda}P(i)/\Tr_{\oplus_{j>i}\,{}_{\Lambda}P(j)}\,({}_{\Lambda}P(i))$. Then,
${}_{\Lambda}\Delta(i)$ is the largest factor module of
${}_{\Lambda}P(i)$ with composition factors only amongst
${}_{\Lambda}S(j)$ for $j\leq i.$ 
 The algebra $\Lambda$ is said to be 
a {\bf standardly stratified algebra}, with respect to the linear
order $\leq$ on the set $[1,n]$, if
$\proj\,(\Lambda)\subseteq\F({}_{\Lambda}{\Delta})$ (see \cite{ADL,Dlab,CPS}). In this case, it is also said that the pair $(\Lambda,\leq)$ 
is a standardly stratified algebra (or an ss-algebra for short).
\

Let $\Lambda$ be an algebra and $\leq$ be a linear order on $[1,n],$ where $n=rk\,K_0(\Lambda).$  By \cite{DR}, it follows that the pair 
$({}_\Lambda\Delta,\leq)$ is always a stratifying system (it is known as the {\bf canonical stratifying system}). Furthermore, if $(\Lambda,\leq)$ is an  ss-algebra, the
representative set of the indecomposable projective $\Lambda$-modules ${}_\Lambda P=\{{}_\Lambda P(i)\}_{i=1}^n$ satisfies that the triple $({}_\Lambda\Delta,{}_\Lambda P,\leq)$ is the Ext-projective 
stratifying system associated to $({}_\Lambda\Delta,\leq).$
\

The main connection between Ext-projective stratifying systems and the class of ss-algebras is given by the following result.

\begin{teo}\label{TeoEpss}\cite{MMS2} \label{teo. sist. estratif. Ext-proyectivos MMS}
 Let $(\Theta,\underline{Q},\leq)$ be an Ext-projective stratifying
system of size $t$ in $\modu\,(\Lambda),$ $\Gamma=\End_\Lambda({Q})^{op}$,
 $H=\Hom_\Lambda(Q,-):\modu\,(\Lambda)\rightarrow \modu\,(\Gamma)$ and
 $L= Q\otimes_\Gamma -:\modu\,(\Gamma)\rightarrow \modu\,(\Lambda)$.
 Then, the following statements hold true.
\begin{itemize}
\item[(a)] The family ${}_\Gamma P=\{H(Q(i)):i\in[1,t]\}$ is a
representative set of the indecomposable projective
$\Gamma$-modules. In particular, $\Gamma$ is a basic algebra and
$rk\,K_0(\Gamma)=t.$
\vspace{.2cm}
\item[(b)] $(\Gamma,\leq)$ is an ss-algebra, that is, $\proj\,(\Gamma)\subseteq
\F({}_{\Gamma}\Delta).$
\vspace{.2cm}
\item[(c)] The restriction $H|_{\F(\Theta)}:\F(\Theta)\to \F({}_{\Gamma}\Delta)$ is an exact equivalence
of categories and
$L|_{\F({}_{\Gamma}\Delta)}:\F({}_{\Gamma}\Delta)\to
\F(\Theta)$ is a quasi-inverse of $H|_{\F(\Theta)}.$
\vspace{.2cm}
\item[(d)] $H(\Theta(i))\simeq{}_{\Gamma}\Delta(i)$, for all $i\in[1,t]$.
\vspace{.2cm}
\item[(e)] $\add\,(Q)=\F(\Theta)\cap {}^{\perp}\F(\Theta),$ where $M\in{}^{\perp}\F(\Theta)$ if and only if the restriction functor
 $\Ext_\Lambda^1(M,-)|_{\F(\Theta)}=0.$
\end{itemize}
\end{teo}
\

Another nice feature, for a stratifying system $(\Theta,\leq)$ of size $t,$ is that an analogous of the 
Jordan-Holder's Theorem holds for the set of ``relative simples'' $\Theta$ in $\F(\Theta).$ That is, for 
any $M\in\F(\Theta)$ and all $i\in[1,t],$ the filtration multiplicity $[M:\Theta(i)]$ is well defined (see 
\cite[Lemma 2.6 (c)]{MMS2}). Therefore, we have the so-called $\Theta$-length $\ell_\Theta(M):=\sum_{i=1}^t\,[M:\Theta(i)]$ of $M.$
\vspace{.3cm}

In what follows, we introduce some features about  split-by-nilpotent extensions.

\begin{defi}\cite{AM,M} Let $\Gamma$ and $\Lambda$ be algebras, and let $I$ be a two-sided ideal of 
$\Gamma.$ It is said that $\Gamma$ is a {\bf split-by-nilpotent extension} of $\Lambda$ by $I,$ 
if $I\subseteq \rad\,(\Gamma)$ and there is an exact sequence of abelian groups 
$$\xymatrix{0\ar[r]& I\ar[r] & \Gamma\ar[r]^\pi & \Lambda\ar[r] & 0}$$
such that $\pi$ is an epimorphism of algebras and there is a morphism of algebras 
$\sigma: \Lambda \rightarrow \Gamma$ with $\pi \sigma = 1_{\Lambda}$.
\end{defi}

In all that follows, we fix the two algebras $\Lambda$ and $\Gamma,$ and the two-sided ideal 
$I\unlhd \Gamma,$ such that $\Gamma$ is a split-by-nilpotent extension of $\Lambda$ 
by $I.$ Observe that the morphisms of algebras $\sigma: \Lambda \rightarrow \Gamma$ and 
$\pi: \Gamma \rightarrow \Lambda,$ induce in a natural way (change of rings), a bimodule structure on $I,$ $\Gamma$
and $\Lambda.$ Furthermore, $\Gamma =\Lambda \oplus I$ as abelian groups, and the multiplicative 
structure of $\Gamma$ can be seen as $\gamma_1 \gamma_2 = (\lambda_1 , i_1 )(\lambda_2 ,
i_2 ) = (\lambda_1 \lambda_2 , i_1 \lambda_2 +\lambda_1 i_2 + i_1 i_2 );$  in this 
case $\pi(\lambda , i )= \lambda$ and $\sigma(\lambda ) = (\lambda , 0).$
\

 We also remark that in \cite{ACT}, the authors consider a quotient path algebra $\Gamma$ and give sufficient conditions 
on a set of arrows $\mathcal{A}$ of the ordinary quiver of $\Gamma,$ so that $\Gamma$ is a split-by-nilpotent extension of 
$\Lambda:=\Gamma/I$ by $I,$  where $I$ is the ideal of $\Gamma$ generated by the set $\mathcal{A}.$ 

\begin{rk}\label{basicos1} The morphisms of algebras $\sigma: \Lambda \rightarrow \Gamma$ and 
$\pi: \Gamma \rightarrow \Lambda$ have the following properties.
\begin{itemize}
\item[(a)]  $\pi$ is a morphism of $\Gamma-\Gamma$ bimodules. In particular, we have the exact sequence of $\Gamma-\Gamma$ bimodules
\[ \xymatrix{0\ar[r]& {}_\Gamma I_\Gamma\ar[r] & {}_\Gamma\Gamma_\Gamma\ar[r]^\pi & {}_\Gamma\Lambda_\Gamma\ar[r] & 0.}\]
\item[(b)] $\pi$ and $\sigma$ are morphisms of $\Lambda-\Lambda$ bimodules. In particular
\[ {}_{\Lambda}\Gamma_{\Lambda}={}_{\Lambda}\Lambda_{\Lambda}\bigoplus{}_{\Lambda}I_{\Lambda}\]
as $\Lambda-\Lambda$ bimodules.
\item[(c)] ${}_{\Lambda}\Gamma_{\Gamma}\otimes_{\Gamma}{}_{\Gamma}\Lambda_{\Lambda}\simeq  {}_{\Lambda}\Lambda_{\Lambda}\simeq {}_\Lambda\Lambda_\Gamma\otimes_{\Gamma}{}_{\Gamma}\Gamma_{\Lambda}$ as  $\Lambda-\Lambda$ bimodules. 
\end{itemize}
\end{rk}

\begin{lem}\label{CambioAni1}
Let $M\in \modu\,(\Gamma)$ and consider ${}_\Lambda M$ as $\Lambda$-module given by the change of rings $\sigma :\Lambda \rightarrow \Gamma$. Then, there exist natural isomorphisms
\[\Hom_{\Gamma}({}_\Gamma \Gamma_{\Lambda} ,{}_{\Gamma}M) \simeq {}_\Lambda M \simeq \Hom_{\Gamma} ({}_{\Gamma} \Lambda_{\Lambda},{}_{\Gamma}M)\]
\end{lem}
\begin{dem} It is straightforward to see that the natural morphisms $$\varphi_M:\Hom_{\Gamma}({}_\Gamma \Gamma_{\Lambda} ,{}_{\Gamma}M)\to{}_{\Lambda}M\;\text{ and  }\;\psi_M:{}_{\Lambda}M\to\Hom_{\Gamma} ({}_{\Gamma} \Lambda_{\Lambda},{}_{\Gamma}M),$$ given by $\varphi_M (f):=f(1)$ y $\psi_M(m)(\lambda)=\sigma (\lambda)m$, are isomorphisms of $\Lambda$-modules such that $\psi_{M}^{-1}(f) =f(1)$ and $\varphi_{M}^{-1} (m)(\gamma )= \gamma m.$
\end{dem}

\begin{lem}\label{CambioAni2} $\Hom_{\Lambda} ({}_{\Lambda}\Lambda_{\Gamma}, {}_{\Lambda}\Lambda_{\Gamma} ) \simeq {}_{\Gamma}\Lambda_{\Gamma}$ as $\Gamma-\Gamma$ bimodules.
\end{lem}
\begin{dem} The morphism $\varphi:\Hom_{\Lambda} ({}_{\Lambda}\Lambda_{\Gamma}, {}_{\Lambda}\Lambda_{\Gamma} ) \rightarrow{}_{\Gamma}\Lambda_{\Gamma}$, given by $\varphi (f) := f(1)$ is an isomorphism of $\Gamma -\Gamma$ bimodules with inverse $\varphi^{-1} (\lambda )(x)= x\lambda.$
\end{dem}

\section{The usual change of rings functors}

Let $\Gamma$ be a split-by-nilpotent extension  of $\Lambda$ by  $I\unlhd\Gamma$. 
We have the functors
\begin{displaymath}
  \xymatrix{\modu\,(\Gamma)\ar[r]^{F} & \modu\,(\Lambda)\ar[r]^{G} & \modu\,(\Gamma),}
\end{displaymath}
where $F:={}_{\Lambda}\Lambda_{\Gamma} \otimes_{\Gamma} -$ and $G:=\ _{\Gamma}\Gamma_{\Lambda} \otimes_{\Lambda} -.$ These functors are known as 
{\bf {change of rings functors}}. 
\

 We also recall, that a functor $H:\mathcal{A}\to \mathcal{B},$ between additive categories, 
{\bf{reflects}} zero objects if $H(A)=0$ implies that $A=0$ for any $A\in\mathcal{A}.$ Furthermore, for a given class $\X$  of objects in 
$\mathcal{A},$ {\bf{the essential image}} of the functor $H|_\X:\X\to\mathcal{B},$ which is denoted by $\Ima\,(H|_\X)$ or $H(\X),$ is the full subcategory of $\mathcal{B}$ whose objects are  all the objects $Z\in\mathcal{B}$ for which
there is an object $X\in\X$ such that $Z\simeq H(X).$ 
\

In the following lemma, we write down some well-known basic properties (see \cite{AM, AZ, MMerkS}), and for the convenience of the reader, it is included a proof.

\begin{lem}\label{basicos2} For the change of rings functors $F$ and $G,$ the following statements hold true.
\begin{itemize}
\item[(a)] $FG\simeq 1_{\modu\,(\Lambda)}.$
\item[(b)] The functors $F$ and $G$ are faithful. In particular, they reflect zero objects. 
\item[(c)] For any $M\in\modu\,(\Lambda),$ the $\Gamma$-module $G(M)$ is indecomposable if and only if $M$ is indecomposable.
\item[(d)] For any $N\in\modu\,(\Gamma),$ if the $\Lambda$-module $F(N)$ is indecomposable then $N$ is indecomposable.
\end{itemize}
\end{lem} 
\begin{dem} (a) It follows from \ref{basicos1} (c).
\

(b) The fact that $G$ is faithful follows from (a). Let us prove that $F$ is also faithful. Indeed, for any $M\in\modu\,(\Gamma),$ it can be seen easily that $\varphi_M:{}_\Lambda\Lambda
_\Gamma\otimes_\Gamma M\to {}_\Lambda M,$ given by $\varphi_M(\lambda\otimes m):=\sigma(\lambda)m,$ is an isomorphism of $\Lambda$-modules, where ${}_\Lambda M$ has the structure of $\Lambda$-module given by the change of rings $\sigma :\Lambda \rightarrow \Gamma.$
Furthermore, for any $f:M\to N$ in $\modu\,(\Gamma),$ we have the following commutative diagram
$$\xymatrix{F(M)\ar[r]^{F(f)}\ar[d]_{\varphi_M} & F(N)\ar[d]^{\varphi_N}\\
{}_\Lambda M\ar[r]_f & {}_\Lambda N.}$$ 
Thus, if $F(f)=0$ then $f=0;$ proving that $F$ is faithful.
\

(c) Let $M\in\modu\,(\Lambda)$ be such that $G(M)$ is indecomposable. In particular, $M\neq 0$. If $M=M_1 \oplus M_2$ then 
$G(M)=G(M_1)\oplus G(M_2);$ and since $G(M)$ is indecomposable, we have that $G(M_1)=0$ or $G(M_2)=0$. Thus, by (b), it follows that $M_1=0$ or 
$M_2=0;$ proving that $M$ is indecomposable.\\
Let $M\in\modu\,(\Lambda)$ be an indecomposable $\Lambda$-module. So $M\neq 0$ and by (b) $G(M)\neq 0$. If $G(M)= N_1 \oplus N_2$ then 
by (a) $M = F(N_1)\oplus F(N_2)$. Therefore, using that $M$ is indecomposable, it follows that  $F(N_1)=0$ or $F(N_2)=0$. Thus, by (b), 
we have that  $N_1=0$ or $N_2=0;$ and so $G(M)$ is indecomposable. 
\

(d) As in (c), the item (d) follows from the fact that $F$ reflects zero objects. 
\end{dem}

\begin{lem}\label{basicoG} Let $I_{\Lambda}$ be a projective $\Lambda$-module. Then, the following statements hold true.
\begin{itemize}
\item[(a)] $\Gamma_\Lambda$ is a projective $\Lambda$-module and $\Tor_1^\Gamma(\Lambda_\Gamma,-)|_{\Ima\,(G)}=0.$
\vspace{.2cm}
\item[(b)] For all $n\geq 0$ and any $X,Y \in \modu\,(\Lambda )$, we have that 
\[ \Ext_{\Gamma}^{n}  (G(X),G(Y))\simeq \Ext_{\Lambda}^{n}  (X,Y) \oplus \Ext_{\Lambda}^{n}  (X,I \otimes_{\Lambda} Y).\]
\end{itemize}
\end{lem}
\begin{dem} (a) Since $\Gamma_\Lambda=\Lambda_\Lambda\oplus I_\Lambda$ (see \ref{basicos1} (b)), it follows that $\Gamma_\Lambda$ is a projective 
$\Lambda$-module, and hence $\Tor_1^\Lambda(\Gamma_\Lambda,-)=0.$ In particular, we get that $G$ is an exact functor.
\ 

Let $X\in\Ima\,(G).$ Then there is an isomorphism $f:G(M)\to X$ for some $M\in\modu\,(\Lambda).$ Consider an exact sequence ending at $M,$ that is, 
$\eta\;:\;0\to K\to P\stackrel{h}{\to} M\to 0.$ Since $G$ is an exact functor, we get the exact sequence 
$\eta'\;:\;0\to G(K)\to G(P)\stackrel{fG(h)}{\to} X\to 0.$ Thus, by applying the functor $F$ to $\eta'$ and using \ref{basicos2} (a), we get an 
exact and commutative diagram 
$$\xymatrix{0\ar[r] & \Tor_1^\Gamma(\Lambda_\Gamma,X)\ar[r]& FG(K)\ar[d]\ar[r] & FG(P)\ar[d]\ar[r] & FG(X)\ar[d]\ar[r] & 0\\
& 0\ar[r] & K \ar[r] & P\ar[r] & M\ar[r] & 0,}$$
where the vertical arrows are isomorphisms; proving that $\Tor_1^\Gamma(\Lambda_\Gamma,X)=0.$
\

(b) Since $\Gamma_\Lambda$ is projective, by \cite[Exercise 9.21]{Ro}, we get that $$\Ext_{\Gamma}^{n}(G(X),G(Y))\simeq \Ext_{\Lambda}^{n}(X,
\Hom_\Gamma({}_\Gamma\Gamma_\Lambda,G(Y))).$$ Thus, the result follows from \ref{basicos1} (b) and \ref{CambioAni1}.
\end{dem}

\begin{teo}\label{EquiFG} For the change of rings functors $F$ and $G,$ the following statements hold true.
\begin{itemize}
\item[(a)] The restriction functor $F|_{\Ima\,(G)}:\Ima\,(G)\to \modu\,(\Lambda)$ is an  equivalence and $G:\modu\,(\Lambda)\to \Ima\,(G)$ is a quasi-inverse. Moreover, if  
$I_\Lambda$ is projective, then  $F|_{\Ima\,(G)}$ and $G$ are exact functors. 
\item[(b)] $\add\,(GM)\subseteq \Ima\,(G)$ for any $M\in\modu\,(\Lambda).$ Thus, the restriction functor $F|_{\add\,(GM)}:\add\,(GM)\to \add\,(M)$ is an equivalence and  
$G|_{\add\,(M)}:\add\,(M)\to \add\,(GM)$ is a quasi-inverse.  
\end{itemize}
\end{teo}
\begin{dem} (a) Let $\varepsilon:FG\to 1_{\modu\,(\Lambda)}$ be the isomorphism of functors given in \ref{basicos2} (a). Firstly, we assert that $F|_{\Ima\,(G)}$ is full. Indeed, let $X,Y\in\Ima\,(G).$ Then there are isomorphisms $\alpha_1:X\to G(M)$ and  $\alpha_2:Y\to G(N)$ 
for some $M,N\in\modu\,(\Lambda).$ We need to show that $F:\Hom_\Gamma(X,Y)\to \Hom_\Lambda(FX,FY)$ is surjective. For any $f\in \Hom_\Lambda(FX,FY),$ set $\overline{f}:=
F(\alpha_2)fF(\alpha_1^{-1}),$ $f':=\varepsilon_N\overline{f}\varepsilon_M^{-1}$ and $h:=
\alpha_2^{-1}G(f')\alpha_1.$ So an straightforward calculation gives us that $F(h)=f,$ proving 
that $F|_{\Ima\,(G)}$ is full. Observe, that $F|_{\Ima\,(G)}$ is dense since $FG\simeq 1_{\modu\,(\Lambda)}.$ Moreover, by \ref{basicos2} (b)  it follows that 
$F|_{\Ima\,(G)}$ is an equivalence. Furthermore, since $FG\simeq 1_{\modu\,(\Lambda)},$ we conclude that $G:\modu\,(\Lambda)\to \Ima\,(G)$ is a quasi-inverse of $F|_{\Ima\,(G)}.$\\ Finally, assume that $I_\Lambda$ is projective. Then, by \ref{basicoG} (a), we get that $F|_{\Ima\,(G)}$ and $G$ are exact functors. 
\

(b) Let $M\in\modu\,(\Lambda)$ and let $X\in \add\,(GM).$ Then, there exists $Z\in\modu\,(\Gamma)$ such that $X\oplus Z=GM^m,$ for some $m.$ Since $G$ is an additive functor, we may assume that $X$ is indecomposable. Thus, 
from the Krull-Remak-Schmidt Theorem and \ref{basicos2} (c), we get that $X\simeq GM'$ for 
some indecomposable direct summand $M'$ of $M$ and thus  $X\in\Ima\,(G).$
\end{dem}

\begin{cor}\label{EquiProj} The restriction functor $F|_{\proj\,(\Gamma)}:\proj\,(\Gamma)\to 
\proj\,(\Lambda)$ is an equivalence and 
$G|_{\proj\,(\Lambda)}:\proj\,(\Lambda)\to \proj\,(\Gamma)$ is a quasi-inverse.
\end{cor}
\begin{dem} It follows from \ref{EquiFG} (b) since $G({}_\Lambda\Lambda)\simeq {}_\Gamma\Gamma$ and $\add\,({}_\Gamma\Gamma)=\proj\,(\Gamma).$
\end{dem}

\section{The functor $G$ and split-by-nilpotent extensions}

We recall that the term algebra means artin $R$-algebra over a commutative artinian ring $R$ and $\Gamma$ is a split-by-nilpotent extension  of $\Lambda$ by $I$. Consider the change of rings functor $G:=\ _{\Gamma}\Gamma_{\Lambda} \otimes_{\Lambda} -:
\modu\,(\Lambda)\to\modu\,(\Gamma).$ Recall that ${}_{\Lambda}\Gamma_{\Lambda}={}_{\Lambda}\Lambda_{\Lambda}\bigoplus{}_{\Lambda}I_{\Lambda}$ as bimodules. Hence $\Gamma\otimes_\Lambda N=(\Lambda\otimes_\Lambda N)\oplus(I\otimes_\Lambda N)$ as $\Lambda$-modules. For the sake of simplicity, we some times write $(M,N)$ instead of $\Hom_\Lambda(M,N).$ 
\

For each pair $M,N\in \modu\,(\Lambda),$ we consider the isomorphism $\delta=\delta_{M,N}$ of $R$-modules 
\begin{displaymath}
\Hom_\Lambda(M,N)\oplus\Hom_\Lambda(M,I\otimes N)\stackrel{\delta}{\longrightarrow} \Hom_\Gamma(GM,GN), 
\end{displaymath}
which is obtained as the composition $\delta:=\tau\circ\Hom_\Lambda(M,\varphi^{-1})\circ\nu$
of the following isomorphism of $R$-modules:
\begin{itemize}
\item[(a)] $\nu:(M,N)\oplus(M,I\otimes N)\to (M,\Lambda\otimes N)\oplus(M,I\otimes N)=(M,\Gamma\otimes N),$ where $\nu(f,g):=(\overline{f},g)$ and $\overline{f}
(m):=1\otimes f(m)$ for all $m\in M.$
\item[(b)] $\Hom_\Lambda(M,\varphi^{-1}):(M,\Gamma\otimes N)\to(M,\Hom_\Gamma(\Gamma,GN)),$ where 
$$\varphi^{-1}(\gamma_1\otimes n)(\gamma_2)=\gamma_2(\gamma_1\otimes n)=\gamma_2\gamma_1\otimes n $$ as can be seen from the proof of \ref{CambioAni1}.
\item[(c)]$\tau:(M,\Hom_\Gamma(\Gamma,GN))\to\Hom_\Gamma(\Gamma\otimes M,GN)=\Hom_\Gamma(GM,GN),$ where $\tau(\alpha)(\gamma\otimes m)=\alpha(m)(\gamma).$
\end{itemize}

\begin{pro}\label{MatrixDelta} $\delta=(G,L):(M,N)\oplus(M,I\otimes N)\to \Hom_\Gamma(GM,GN)$ as a matrix, where 
$L(g)(\gamma\otimes m)=\gamma g(m).$ Furthermore $L:\Hom_\Lambda(M,I\otimes N)\to \Hom_\Gamma(GM,GN)$ is a monomorphism.
\end{pro}
\begin{dem} Consider the natural inclusion $i_1:(M,N)\to (M,N)\oplus(M,I\otimes N).$ We assert that $\delta i_1=G.$ Indeed, let $f\in(M,N),$ $\gamma\in\Gamma$ and $m\in M.$ So we have 
$\delta i_1(f)(\gamma\otimes m)=\tau(\varphi^{-1}\overline{f})(\gamma\otimes m)=(\varphi^{-1}\overline{f}(m))(\gamma)=\gamma(1\otimes f(m))=\gamma\otimes f(m)=G(f)(\gamma\otimes m).$
\

Let $i_2:(M,I\otimes N)\to (M,N)\oplus(M,I\otimes N)$ be the natural inclusion. We check that 
$\delta i_2=L.$ Indeed, let $g\in(M,I\otimes N),$ $\gamma\in\Gamma$ and $m\in M.$ So we have
$\delta i_2(g)(\gamma\otimes m)=(\tau\varphi^{-1}(0,g))(\gamma\otimes m)=\varphi^{-1}(0,g)(m)(
\gamma)=\gamma g(m).$
\

Finally, let $y\in \Hom_\Lambda(M,I\otimes N)$ be such that $L(y)=0.$ Then $\delta\begin{pmatrix}0\\y \end{pmatrix}=G(0)+L(y)=0$ and since 
$\delta$ is an isomorphism, it follows that $y=0 ;$ proving that $L$ is a monomorphism.
\end{dem}
\begin{cor}\label{FielPlenoG} Let $M,N\in \modu\,(\Lambda).$ Then, the following conditions 
are equivalent.
\begin{itemize}
\item[(a)] $G:\Hom_\Lambda(M,N)\to\Hom_\Gamma(GM,GN)$ is an isomorphism.
\item[(b)] $\Hom_\Lambda(M,I\otimes N)=0.$
\end{itemize}  
\end{cor}
\begin{dem} Consider the natural inclusion $i_1:(M,N)\to (M,N)\oplus(M,I\otimes N).$ Then, by \ref{MatrixDelta}, we know that $\delta i_1=G$ and thus the result follows since $\delta$ is an 
isomorphism. 
\end{dem}

\begin{teo}\label{SplitAlgG} For any $M\in\modu\,(\Lambda),$ the algebra $\End_\Gamma(GM)$ is an split-by-nilpotent extension of 
$\End_\Lambda(M)$ by $\Hom_\Lambda(M,I\otimes M).$
\end{teo}
\begin{dem} Let $M\in\modu\,(\Lambda).$ By \ref{MatrixDelta}, we have the exact sequence of $R$-modules 
$$0\longrightarrow\Hom_\Lambda(M,I\otimes M)\stackrel{L}{\longrightarrow}\End_\Gamma(GM)\stackrel{\vartheta}{\longrightarrow}\End_\Lambda(M)
\longrightarrow 0,$$ where $\vartheta:=\pi_1\delta^{-1}$ and $\pi_1:\End_\Lambda(M)\oplus\Hom_\Lambda(M,I\otimes M)\to \End_\Lambda(M)$ is the 
canonical projection. Observe that $G:\End_\Lambda(M)\to\End_\Gamma(GM)$ is a ring morphism and $\vartheta G$ is the identity map.
\

We assert that $\vartheta:\End_\Gamma(GM)\to\End_\Lambda(M)$ is a ring morphism. Indeed, we firstly transfer, by using the isomorphism $\delta,$ 
the multiplicative structure of the ring $\End_\Gamma(GM)$ to the $R$-module $\End_\Lambda(M)\oplus\Hom_\Lambda(M,I\otimes M).$ That is, for any 
$\alpha=(f_\alpha,g_\alpha)$ and $\beta=(f_\beta,g_\beta)$ in $\End_\Lambda(M)\oplus\Hom_\Lambda(M,I\otimes M),$ we set $\alpha\beta:=
\delta^{-1}(\delta(\alpha)\delta(\beta)).$ In what follows, we shall compute the above product and show that 
$$(*)\quad\quad\delta(\alpha)\delta(\beta)=G(f_\alpha f_\beta)+L(g_\alpha f_\beta+\varepsilon).$$ 
Indeed $\delta(\alpha)\delta(\beta)=G(f_\alpha)G(f_\beta)+L(g_\alpha)G(f_\beta)+L(g_\alpha)L(g_\beta)+G(f_\alpha)
L(g_\beta).$  But $L(g_\alpha)G(f_\beta)=L(g_\alpha f_\beta)$ since $L(g_\alpha)G(f_\beta)(\gamma\otimes m)=L(g_\alpha)(\gamma\otimes 
f_\beta(m))=\gamma g_\alpha(f_\beta(m))=L(g_\alpha f_\beta)(\gamma\otimes m).$ To compute $\mu:=L(g_\alpha)L(g_\beta)+G(f_\alpha)
L(g_\beta),$ we proceed as follows. Observe that $\mu=(G(f_\alpha)+L(g_\alpha))L(g_\beta)=\delta(\alpha)L(g_\beta).$ Consider the morphism 
$(0,g_\beta):M\to GM.$ Using the fact that $I\unlhd\Gamma,$ it can be seen that $\Ima\,(\delta(\alpha)(0,g_\beta))\subseteq I\otimes M.$ Thus, the 
morphism $\varepsilon,$ which is the composition of $M\to \Ima\,(\delta(\alpha)(0,g_\beta))\subseteq I\otimes M,$ satisfies that 
$L(\varepsilon)=\mu.$ Indeed, $L(\varepsilon)(\gamma\otimes m)=\gamma\varepsilon(m)=\gamma\delta(\alpha)(0,g_\beta)(m)=\gamma\delta(\alpha)(
g_\beta(m))=\delta(\alpha)(\gamma g_\beta(m))=\delta(\alpha)(L(g_\beta)(\gamma\otimes m))=\delta(\alpha)L(g_\beta)(\gamma\otimes m);$ proving 
$(*).$ Now, we are 
ready to prove that $\vartheta:\End_\Gamma(GM)\to\End_\Lambda(M)$ is a ring homomorphism. That is, by $(*),$ we have that $\vartheta(\delta(\alpha)
\delta(\beta))=\pi_1\delta^{-1}(\delta(\alpha)\delta(\beta))=f_\alpha f_\beta=\pi_1\delta^{-1}\delta(\alpha)\pi_1\delta^{-1}\delta(\beta)=
\vartheta(\delta(\alpha))\vartheta(\delta(\beta)).$
\

Finally, we prove that $\Ima\,(L)\subseteq\rad\,(\End_\Gamma(GM)).$ To see that, it is enough to check that the ideal $\Ima\,(L)$ is nilpotent. 
Let $g_1,g_2,\cdots,g_n\in\Hom_\Lambda(M,I\otimes M),$ $\gamma\in\Gamma$ and $m\in M.$ Since $L(g_1)L(g_2)\cdots L(g_n)(\gamma\otimes m)
\in I^n\otimes M;$ and using the fact that $I$ is nilpotent, it follows that $\Ima\,(L)$ is also nilpotent.
\end{dem}

\section{Extending stratifying systems with the functor $G$}

In this section, we consider a split-by-nilpotent extension $\Gamma$ of $\Lambda$ by $I$. As we 
have seen before, there is the change of rings functor $G:=\ _{\Gamma}\Gamma_{\Lambda} \otimes_{\Lambda} -:
\modu\,(\Lambda)\to\modu\,(\Gamma).$ We give conditions for the image under $G,$ of a stratifying system in $\modu\,(\Lambda),$ to be a stratifying system in $\modu\,(\Gamma).$

\begin{defi}\label{DefCompatible} A stratifying system $(\Theta,\leq),$ of size $t$ in $\modu\,(\Lambda),$ is {\bf 
compatible} with the ideal $I\unlhd \Gamma$ if the following conditions hold.
\begin{itemize}
\item[(C1)] $\Hom_\Lambda(\Theta(j),I\otimes_{\Lambda}\Theta(i))=0$ for $j>i.$ 
\item[(C2)] $\Ext_\Lambda^1(\Theta(j),I\otimes_{\Lambda}\Theta(i))=0$ for $j\geq i.$
\end{itemize}
\end{defi}

\begin{pro}\label{TeoGLevantSS} Let $\Theta =\{ \Theta (i) \}_{i=1}^{t}$ be objects in 
$\modu\,(\Lambda)$ and $\leq$ be a linear order on $[1,t]$. If $I_{\Lambda}$ is projective, then the following conditions are equivalent. 
\begin{enumerate}
\item[(a)] $(G(\Theta), \leq)$ is a stratifying system in $\modu\,(\Gamma)$.
\item[(b)] $(\Theta, \leq)$ is a stratifying system in $mod(\Lambda),$ which is compatible with the ideal $I.$
\end{enumerate}
\end{pro}
\begin{dem} By \ref{basicos2}, we know that the functor 
$G: \modu\,(\Lambda) \rightarrow \modu\,(\Gamma)$ reflects and preserves indecomposables. On 
the other hand, since $I_{\Lambda}$ is projective, we have by \ref{basicoG} that
\[\Ext_{\Gamma}^{i}  (G(X),G(Y)) \simeq \Ext_{\Lambda}^{i}  (X,Y) \oplus \Ext_{\Lambda}^{i}(X, I \otimes_{\Lambda} Y) ,\]
for all $X,Y \in \modu\,(\Lambda)$ and any $i.$ Thus, the equivalence between (a) and (b) follows.
\end{dem}

\begin{cor}\label{CorGLevantSS} Let $\Theta =\{ \Theta (i) \}_{i=1}^{t}$ be objects in 
$\modu\,(\Lambda),$ and let $\leq$ be a linear order on $[1,t]$. If ${}_\Lambda I_{\Lambda}
\in\add\,({}_\Lambda \Lambda_{\Lambda})$ then the following conditions are equivalent. 
\begin{enumerate}
\item[(a)] $(G(\Theta), \leq)$ is a stratifying system in $\modu\,(\Gamma)$.
\item[(b)] $(\Theta, \leq)$ is a stratifying system in $\modu\,(\Lambda).$ 
\end{enumerate}
\end{cor}
\begin{dem} Let ${}_\Lambda I_{\Lambda}\in\add\,({}_\Lambda \Lambda_{\Lambda}).$ Then 
$I \otimes_{\Lambda} X\in\add\,(X)$ for any $X\in\modu\,(\Lambda).$ Therefore, any stratifying 
system in $\modu\,(\Lambda)$ is compatible with the ideal $I.$ So, the result follows from \ref{TeoGLevantSS}.
\end{dem}

Let us consider the following examples.

\begin{ex} Consider the trivial extension $\Gamma:=\Lambda\ltimes I$ of an algebra  $\Lambda$ by its minimal injective cogenerator $I:=D(\Lambda).$ If $\Lambda$ is a symmetric algebra, it is well known that ${}_\Lambda I_\Lambda\simeq {}_\Lambda \Lambda_\Lambda$ as bimodules. Therefore, in this case, the needed hypothesis in \ref{CorGLevantSS} holds. We recall that in \cite{ES} stratifying systems for symmetric special biserial algebras are constructed.
\end{ex}

\begin{ex}\label{Ej2} Let $Q$ be the quiver $\bullet^3\stackrel{\beta}{\to}\bullet^1\stackrel{\alpha}{\to}\bullet^2.$ Consider the quotient path $k$-algebra $\Gamma:=kQ/\la\alpha\beta \ra$ with the ideal $I:=\la\overline{\beta}\ra\unlhd \Gamma.$ Then $\Gamma$ is a split-by-nilpotent extension of $\Lambda:=\Gamma/I$ by $I.$ Moreover, the ordinary quiver $Q_\Lambda$ of $\Lambda$ is $\bullet^3\quad\bullet^1\stackrel{\alpha}{\to}\bullet^2.$ We consider the natural order $1\leq 2\leq 3.$ Then, we have the canonical stratifying system $({}_\Lambda\Delta,\leq)$ in $\modu\,(\Lambda)$ where ${}_\Lambda\Delta(i)={}_\Lambda S(i)$ is the simple $\Lambda$-module associated to the vertex $i\in Q_\Lambda.$
It can be seen that $I_\Lambda\simeq e_3\Lambda$ (where $e_3$ is the primitive idempotent associated with the vertex $3$), $I\otimes_\Lambda{}_\Lambda\Delta(1)=0=I\otimes_\Lambda{}_\Lambda\Delta(2)$ and $I\otimes_\Lambda{}_\Lambda\Delta(3)\simeq{}_\Lambda\Delta(1).$ Thus, the stratifying system 
$({}_\Lambda\Delta,\leq)$ is compatible with the ideal $I;$ and so by \ref{TeoGLevantSS} it follows that $(G({}_\Lambda\Delta), \leq)$ is a stratifying system in $\modu\,(\Gamma).$    
\end{ex}

As we have seen in \ref{TeoGLevantSS}, the notion of stratifying system compatible with the 
ideal $I$ plays an important role. In the following proposition, we give conditions for the canonical stratifying system to be compatible with the ideal $I$. For doing so, let ${}_\Lambda P=\{{}_\Lambda P(i)\}_{i=1}^n$ be a representative set of the indecomposable projective $\Lambda$-modules, where $n:=rk\,K_0(\Lambda),$ and let $\leq$ be a linear order on the set $[1,n]:=\{1,2,\cdots, n\}.$ Let us consider the set of standard $\Lambda$-modules ${}_\Lambda\Delta,$  
computed by using the pair $({}_\Lambda P,\leq ),$ and also     
the representative set ${}_\Lambda S=\{{}_{\Lambda}S(i)\}_{i=1}^n$ of simple $\Lambda$-modules, where  ${}_{\Lambda}S(i):={}_{\Lambda}P(i)/\rad\,({}_{\Lambda}P(i))$ for each $i.$  Recall that each ${}_\Lambda\Delta(i)$ has composition factors only amongst ${}_\Lambda S(j)$ with $j\leq i.$ 
That is, the multiplicity $[{}_\Lambda\Delta(i):{}_\Lambda S(j)]$ of the simple ${}_\Lambda S(j)$ in ${}_\Lambda\Delta(i)$ is equal to zero for $j>i.$ So, we start with the following definition.

\begin{defi} The pair $({}_\Lambda S,\leq)$ is admissible with the ideal $I\unlhd \Gamma$ if 
$$[I\otimes{}_\Lambda S(i):{}_\Lambda S(j)]=0\;\;\;\text{ for }\;\;\;j>i.$$
\end{defi}

\begin{lem}\label{AuxCond2} Let $I_{\Lambda}$ be a projective $\Lambda$-module. Then, the following statements are equivalent.
\begin{itemize}
\item[(a)] The pair $({}_\Lambda S,\leq)$ is admissible with the ideal $I.$
\item[(b)] $[I\otimes{}_\Lambda\Delta(i):{}_\Lambda S(j)]=0$ for $j>i.$
\end{itemize}
\end{lem}
\begin{dem} (a) $\Rightarrow$ (b) It follows from the fact that $I\otimes_\Lambda -:\modu\,(\Lambda)\to \modu\,(\Lambda)$ is an exact 
functor and ${}_\Lambda\Delta(i)\in\F(\{{}_\Lambda S(t)\;:\;t\leq i\})$ for each $i.$
\

(b) $\Rightarrow$ (a) Let $j>i.$ By Applying the exact functor $I\otimes_\Lambda -$ to the  
exact sequence $0\to \rad\,({}_\Lambda\Delta(i))\to {}_\Lambda\Delta(i)\to {}_\Lambda S(i)\to 0,$ 
we get the exact sequence $0\to I\otimes\rad\,({}_\Lambda\Delta(i))\to I\otimes{}_\Lambda\Delta(i)\to I\otimes{}_\Lambda S(i)\to 0.$ Thus, the condition given in (b) implies that $[I\otimes{}_\Lambda S(i):{}_\Lambda S(j)]=0;$ proving (a).
\end{dem}
\vspace{.2cm}

The following result relates the admissibility of $({}_\Lambda S,\leq)$ with the compatibility of 
$({}_\Lambda\Delta,\leq).$

\begin{pro}\label{Condition2} Let $I_{\Lambda}$ be a projective $\Lambda$-module. Then, 
the following statements hold true.
\begin{itemize}
\item[(a)] If $({}_\Lambda S,\leq)$ is admissible with the ideal $I,$ then the canonical stratifying system $({}_\Lambda\Delta,\leq)$ is compatible with the ideal $I.$
\item[(b)] If $\Lambda$ is an ss-algebra such that $({}_\Lambda\Delta,\leq)$ is compatible with the ideal $I,$ then $({}_\Lambda S,\leq)$ is admissible with the ideal $I.$
\end{itemize}
\end{pro}
\begin{dem} (a) $\Rightarrow$ (b) Assume that, for each $i\in[1,n],$ the $\Lambda$-module $I\otimes{}_\Lambda S(i)$ has composition factors only amongst ${}_\Lambda S(j)$ with $j\leq i.$ Then, by \ref{AuxCond2}, we get that $[I\otimes{}_\Lambda\Delta(j):{}_\Lambda S(i)]=0$ for $i>j.$ Therefore $\Hom_\Lambda({}_\Lambda P(i),I\otimes{}_\Lambda\Delta(j))=0$ for $i>j,$ and so the condition (C1) in \ref{DefCompatible} holds.
\

Let now $i\geq j,$ and let $\nu\;:\;0\to I\otimes{}_\Lambda\Delta(j)\to U\stackrel{\alpha}{\to}
{}_\Lambda\Delta(i)\to 0$ be an exact sequence. So, by \ref{AuxCond2} (b), we get that 
$U\in\F(\{{}_\Lambda S(t)\;:\;t\leq i\}).$ Consider the epimorphism $p:{}_\Lambda P(i)\to
{}_\Lambda\Delta(i)$ where $\Ker\,(p)=\Tr_{\oplus_{r>i}\;{}_\Lambda P(r)}\,({}_\Lambda P(i)).$ Then, there is a morphism $f:{}_\Lambda P(i)\to U$ such that $p=\alpha f.$ 
By taking the factorization ${}_\Lambda P(i)\stackrel{\overline{f}}{\to}\Ima\,(f)\stackrel{\imath}
{\to} U$ of $f$ throughout its image, we have that $(\alpha\imath)\overline{f}=p.$ That is, the 
quotient morphism $\overline{f}:{}_\Lambda P(i)\to\Ima\,(f)$ factors throughout 
$p:{}_\Lambda P(i)\to{}_\Lambda\Delta(i).$ Moreover, since $\Ker\,(p)=\Tr_{\oplus_{r>i}\;{}_\Lambda P(r)}\,({}_\Lambda P(i)),$ it follows that  $p$ factors throughout $\overline{f}.$ Hence 
$\alpha\imath:\Ima\,(f)\to {}_\Lambda\Delta(i)$ is an isomorphism. Therefore, the exact sequence 
$\nu$ splits, and so the condition (C2) in \ref{DefCompatible} holds. 
\

(b) $\Rightarrow$ (a) Assume that $\Lambda$ is an ss-algebra such that $({}_\Lambda\Delta,\leq)$ is compatible with the ideal $I.$ Let $j>i$ and consider the canonical exact sequence 
$\eta\;:\;0\to K(j)\to {}_\Lambda P(j)\to {}_\Lambda\Delta(j)\to 0,$ where $K(j):=\Tr_{\oplus_{r>j}\;{}_\Lambda P(r)}\,({}_\Lambda P(j)).$ Since $\Lambda$ is an ss-algebra, it is known that 
$K(j)\in\F(\{{}_\Lambda\Delta(t)\;:\;t>j\}).$ Thus $\Hom_\Lambda(K(j),I\otimes {}_\Lambda\Delta(i))=0$ (see \ref{DefCompatible} (C1)). Applying $\Hom_\Lambda(-,I\otimes {}_\Lambda\Delta(i))$ to $\eta,$ and since $({}_\Lambda\Delta,\leq)$ is compatible with the ideal $I,$ we conclude that $\Hom_\Lambda({}_\Lambda P(j),I\otimes {}_\Lambda\Delta(i))\simeq 
\Hom_\Lambda( K(j),I\otimes {}_\Lambda\Delta(i))=0.$ Therefore $[I\otimes {}_\Lambda\Delta(i):{}_\Lambda S(j)]=0$ for $j>i.$ Finally, by \ref{AuxCond2} (b), we conclude that  $({}_\Lambda S,\leq)$ is admissible with the ideal $I,$
\end{dem}
\vspace{.2cm}

Let $I_{\Lambda}$ be a projective $\Lambda$-module. In \cite{MMerkS}, the authors consider as a main hypothesis that the $\Lambda$-module $I\otimes{}_\Lambda S(i)$ has composition factors only amongst ${}_\Lambda S(j)$ with $j\leq i.$ As we have seen in \ref{Condition2}, under the hypothesis that $\Lambda$ is an ss-algebra, this is equivalent to the compatibility condition (see \ref{DefCompatible}) for the canonical stratifying system. Observe that \ref{DefCompatible} is precisely the needed condition 
to determine when a stratifying system in $\modu\,(\Lambda)$ can be extended, throughout the 
functor $G,$ to a stratifying system in $\modu\,(\Gamma)$  (see \ref{TeoGLevantSS}). 

\begin{pro}\label{FiltrFG} Let $I_{\Lambda}$ be a projective $\Lambda$-module and let 
$(\Theta,\leq)$ be a stratifying system of size $t$ in $\modu\,(\Lambda),$ which 
is compatible with the ideal $I.$ Then, $(G\Theta,\leq)$ is a stratifying system of size $t$ in $\modu\,(\Gamma),$ and the following statements hold true.
\begin{itemize}
\item[(a)] The restriction $F|_{\F(G\Theta)}:\F(G\Theta)\to \F(\Theta)$ is well defined, and it is an exact, faithful and dense functor which reflects indecomposables.
\item[(b)] The restriction $G|_{\F(\Theta)}:\F(\Theta)\to \F(G\Theta)$ is well defined, and it is an exact and faithful functor which preserves and reflects indecomposables. 
\item[(c)] The restriction $F|_{\Ima\,(G|_{\F(\Theta)})}:\Ima\,(G|_{\F(\Theta)})\to \F(\Theta)$ 
is an equivalence of categories, and a quasi-inverse is the restriction $G|_{\F(\Theta)}:\F(\Theta)\to \Ima\,(G|_{\F(\Theta)}).$
\end{itemize}
\end{pro}
\begin{dem} By \ref{TeoGLevantSS} we know that $(G\Theta,\leq)$ is a stratifying in $\modu\,(\Gamma).$
\

(a) Let $M\in\F(G\Theta).$ We prove, by induction on the $G\Theta$-length $\ell_{G\Theta}(M),$ that $F(M)\in\F(\Theta)$ and $\Tor_1^{\Gamma}(\Lambda_\Gamma,M)=0.$ If $\ell_{G\Theta}(M)=1$ then 
$M\simeq G\Theta(i)$ for some $i.$ Thus $F(M)\simeq FG\Theta(i)\simeq \Theta(i)$ (see \ref{basicos2} (a)) and $\Tor_1^{\Gamma}(\Lambda_\Gamma,M)\simeq\Tor_1^{\Gamma}(\Lambda_\Gamma,G\Theta(i))=0$ (see \ref{basicoG} (a)).\\
Let $\ell_{G\Theta}(M)>1.$ Then, from \cite[Lemma 2.8]{MMS2}, there is an exact sequence 
$\eta\;:\;0\to G\Theta(i)\to M\to N\to 0$ in $\F(G\Theta),$ with $\ell_{G\Theta}(N)=\ell_{G\Theta}(M)-1.$ Applying the functor $F$ to $\eta,$ we get the exact sequence 
$\Tor_1^{\Gamma}(\Lambda_\Gamma,G\Theta(i))\to\Tor_1^{\Gamma}(\Lambda_\Gamma,M)\to \Tor_1^{\Gamma}(\Lambda_\Gamma,N)\to FG\Theta(i)\to FM\to FN\to 0.$ By induction we know that 
$\Tor_1^{\Gamma}(\Lambda_\Gamma,N)=0$ and $FN\in\F(\Theta).$ Thus $F(M)\in\F(\Theta)$ and $\Tor_1^{\Gamma}(\Lambda_\Gamma,M)=0.$ In particular, it follows that the restriction $F|_{\F(G\Theta)}$ is well defined and it is also an exact functor. Moreover, by \ref{basicos2}, it 
is also a faithful and dense functor which reflects indecomposables.
\

(b) Since $\Gamma_{\Lambda}$ is projective (see \ref{basicoG} (a)), it follows that $G:\modu\,(\Lambda)\to \modu\,(\Gamma)$ is an exact functor. Hence the restriction $G|_{\F(\Theta)}:\F(\Theta)\to \F(G\Theta)$ 
is well defined. Finally, from \ref{basicos2} (b) and (c), we conclude that $G$ is faithful and also preserves and reflects indecomposables.
\

(c) It follows from (a), (b) and \ref{EquiFG} (a).
\end{dem} 

\begin{cor}\label{Equivalencias} Let $I_{\Lambda}$ be a projective $\Lambda$-module and let 
$(\Theta,\leq)$ be a stratifying system of size $t$ in $\modu\,(\Lambda),$ which 
is compatible with the ideal $I.$ Then, the following conditions are equivalent.
\begin{itemize}
\item[(a)] The restriction functor $F|_{\F(G\Theta)}:\F(G\Theta)\to \F(\Theta)$ is an exact equivalence of categories, and its quasi-inverse is the restriction functor $G|_{\F(\Theta)}:\F(\Theta)\to \F(G\Theta).$
\item[(b)] $\Ima\,(G|_{\F(\Theta)})=\F(G\Theta).$ 
\item[(c)] The class $\Ima\,(G|_{\F(\Theta)})$ is closed under extensions in $\modu\,(\Gamma).$
\end{itemize}
\end{cor}
\begin{dem} It follows from \ref{FiltrFG} and the fact that $\F(G\Theta)$ is the smaller full subcategory of $\modu\,(\Gamma),$ which contains $G\Theta$ and is closed under extensions.
\end{dem}

\vspace{.3cm}

Observe that, in general, the class $\Ima\,(G|_{\F(\Theta)})$ is not necessarily closed under extensions in $\modu\,(\Gamma).$ A sufficient condition for the equality 
$\Ima\,(G|_{\F(\Theta)})=\F(G\Theta)$ will be given in \ref{Gequiv}.

\begin{pro}\label{G-epss} Let $I_{\Lambda}$ be projective, $(\Theta,\leq)$ be a stratifying system of size $t$ in $\modu\,(\Lambda),$ which 
is compatible with the ideal $I,$ and let 
$(\Theta,\underline{Q},\leq)$ be the Ext-projective stratifying system associated to $(\Theta,\leq).$ Then\\
$\Ext_\Lambda^1(Q,I\otimes \mathbf{\Theta})=0$ if and only if the triple $(G\Theta,G\underline{Q},\leq)$ is the Ext-projective stratifying system associated to the stratifying system $(G\Theta,\leq).$
\end{pro}
\begin{dem}  We assert that $\Ext_\Gamma^1(GQ(i),G\Theta(j))\simeq \Ext_\Lambda^1(Q(i),I\otimes\Theta(j))$ for any $i,j\in[1,t].$ Indeed, since $\Ext_\Lambda^1(Q(i),\Theta(j))=0$ for any $i,j,$ then by \ref{basicoG} (b), the assertion follows. 
Thus, by the above assertion, the implication "$\Leftarrow$" is clear. Assuming that $\Ext_\Lambda^1(Q,I\otimes \mathbf{\Theta})=0,$ we obtain from our assertion that $GQ$ is Ext-projective in $G\Theta.$ Furthermore, using the fact 
that $G|_{\F(\Theta)}$ is an exact functor (see \ref{FiltrFG} (b)), we get that the fundamental sequence $\varepsilon_i\;:\;0\to K(i)\to Q(i)\to \Theta(i)\to 0,$ attached to the system 
$(\Theta,\underline{Q},\leq)$ (see \ref{definicion de epps} (b)) gives the fundamental sequence $G\varepsilon_i\;:\;0\to GK(i)\to GQ(i)\to G\Theta(i)\to 0$ corresponding to the 
system $(G\Theta,G\underline{Q},\leq).$
\end{dem}

\begin{teo}\label{GbarraInducido} Let $I_{\Lambda}$ be projective, $(\Theta,\leq)$ be a stratifying system of size $t$ in $\modu\,(\Lambda),$ 
which is compatible with the ideal $I$ and $\Ext_\Lambda^1(Q,I\otimes \mathbf{\Theta})=0,$ 
and let $(\Theta,\underline{Q},\leq)$ be the Ext-projective stratifying system associated to $(\Theta,\leq).$ Consider the algebras 
$A:=\End_\Lambda(Q)^{op}$ and $GA:=\End_\Gamma(GQ)^{op},$ where $Q:=\oplus_{i=1}^t\,Q(i).$ Then, the following statements hold true.
\begin{itemize}
 \item[(a)]  $GA$ is an split-by-nilpotent extension of $A$ by $\Hom_\Lambda(Q,I\otimes Q).$ Furthermore, both algebras $A$ 
and $G(A)$ are basic, standardly stratified and $rk\,K_0(A)=t=rk\,K_0(GA).$
\item[(b)] $\overline{G}:=\Hom_\Gamma(GQ,-)\circ G|_{\F(\Theta)}\circ Q\otimes_A -:\F({}_A\Delta)\to \F({}_{GA}\Delta)$ is well defined and it is an exact and faithful functor which preserves and reflects indecomposables, 
and $\overline{G}({}_A\Delta(i))\simeq {}_{GA}\Delta(i)$ for any $i\in[1,t].$
\item[(c)] $\overline{F}:=\Hom_\Lambda(Q,-)\circ F|_{\F(G\Theta)}\circ GQ\otimes_{GA} -:\F({}_{GA}\Delta)\to \F({}_{A}\Delta)$ is well defined and it is an exact, dense and faithful functor which reflects indecomposables, 
and $\overline{F}({}_{GA}\Delta(i))\simeq {}_{A}\Delta(i)$ for any $i\in[1,t].$
\item[(d)] $A$ is quasi-hereditary if and only if $GA$ is quasi-hereditary.
\item[(e)] The restriction functor $\overline{F}|_{\proj\,(GA)}:\proj\,(GA)\to \proj\,(A)$ is an 
equivalence and its quasi-inverse is $\overline{G}|_{\proj\,(A)}:\proj\,(A)\to \proj\,(GA).$
\end{itemize}
\end{teo}
\begin{dem} (a) It follows from \ref{SplitAlgG} and \cite[Theorem 3.2 (a)]{MMS2}.
\ 

(b) By \cite[Theorem 3.2 ]{MMS2}, it follows that the functors $Q\otimes_A -:\F({}_A\Delta)\to \F(\Theta)$ and $\Hom_\Gamma(GQ,-):\F(G\Theta)\to
\F({}_{GA}\Delta)$ are exact equivalences. Thus, by \ref{FiltrFG} (b), the functor 
 $\overline{G}:=\Hom_\Gamma(GQ,-)\circ G|_{\F(\Theta)}\circ Q\otimes_A -:\F({}_A\Delta)\to \F({}_{GA}\Delta)$ is well defined and has the desired properties. Moreover, by  \cite[Theorem 3.1 ]{MMS2}, we have that $Q_A\otimes{}_A\Delta(i)\simeq \Theta(i)$ and 
$\Hom_\Gamma(GQ,GQ(i))\simeq {}_{GA}\Delta(i).$ Therefore $\overline{G}({}_A\Delta(i))\simeq {}_{GA}\Delta(i)$ for any $i\in[1,t].$
\

(c) It follows, as in the proof of (b), from \cite[Theorem 3.2 ]{MMS2} and \ref{FiltrFG} (a).
\

(d) By \cite[Theorem 3.2]{MMS2}, it is enough to see that: $\rad\,(\End_\Lambda(\Theta(i)))=0$ 
if and only if $\rad\,(\End_\Gamma(G\Theta(i)))=0$ for any $i\in[1,t].$ By \ref{FiltrFG} (a) and (b), we have the ring morphisms $\End_\Lambda(\Theta(i))\stackrel{G}{\to}\End_\Gamma(G\Theta(i))
\stackrel{F}{\to}\End_\Lambda(\Theta(i)).$ Thus $G(\rad\,(\End_\Lambda(\Theta(i))))\subseteq 
\rad\,(\End_\Gamma(G\Theta(i)))$ and also we have the inclusion\\ $F(\rad\,(\End_\Gamma(G\Theta(i))))\subseteq 
\rad\,(\End_\Lambda(\Theta(i))).$ Therefore, the desired equivalence holds by the fact that $F$ and $G$ are faithful functors.
\

(e) It follows from \ref{EquiFG} (b) and \cite[Theorem 3.2]{MMS2}. 
\end{dem}

\begin{rk}\label{FilItT} Let $(\Theta,\underline{Q},\leq)$ be an Ext-projective stratifying 
system in $\modu\,(\Lambda).$ Observe that, if  $I\otimes \mathbf{\Theta}\in\F(\Theta)$ 
then $\Ext_\Lambda^1(Q,I\otimes \mathbf{\Theta})=0.$
\end{rk}

As we have seen above (see \ref{FiltrFG} and \ref{Equivalencias}) the restriction functor $G|_{\F(\Theta)}:\F(\Theta)\to \F(G\Theta)$ is not, in general, 
an equivalence. So, in the following results, we give sufficient conditions ensuring that $G|_{\F(\Theta)}$ is an equivalence of categories. 
In order to do that, we start with the following lemma.

\begin{lem}\label{GfielPleno} Let  $(\Theta,\underline{Q},\leq)$ be an Ext-projective stratifying system of size $t,$ in 
$\modu\,(\Lambda),$ such that $\Hom_\Lambda(\mathbf{\Theta},I\otimes\mathbf{\Theta})=0.$ Then, the following statements hold true.
\begin{itemize}
 \item[(a)]  $\Hom_\Lambda(M,I\otimes N)=0$ for any $M,N\in\F(\Theta).$
 \item[(b)]  $\Ext_\Lambda^1(Q,I\otimes\mathbf{\Theta})=0\quad\Leftrightarrow\quad \Ext_\Lambda^1(\mathbf{\Theta},I\otimes \mathbf{\Theta})=0.$
\end{itemize}
\end{lem}
\begin{dem} (a) It follows from \cite[Lemma 2.8]{MMS2} and by induction on the $\Theta$-length $\ell_\Theta(M)$ of $N.$ 
\

(b) The implication ``$\Leftarrow$'' follows easily since $Q\in\F(\Theta).$ Let $\Ext_\Lambda^1(Q,I\otimes \mathbf{\Theta})=0.$ Then, by 
applying $\Hom_\Lambda(-,I\otimes\Theta(j))$ to the canonical exact sequence $0\to K(i)\to Q(i)\to \Theta(i)\to 0$ in $\F(\Theta),$ we get that 
$\Ext_\Lambda^1(\Theta(i),I\otimes\Theta(j))=0$ since  $\Hom_\Lambda(K(i),I\otimes\Theta(j))=0=\Ext_\Lambda^1(Q(i),I\otimes\Theta(j));$ and 
hence $\Ext_\Lambda^1(\mathbf{\Theta},I\otimes \mathbf{\Theta})=0.$ 
\end{dem}

\begin{teo}\label{Gequiv}  Let $I_{\Lambda}$ be projective, $(\Theta,\underline{Q},\leq)$ be an Ext-projective stratifying system of size $t,$ in 
$\modu\,(\Lambda),$ such that $\Hom_\Lambda(\mathbf{\Theta},I\otimes\mathbf{\Theta})=0=\Ext_\Lambda^1(Q,I\otimes\mathbf{\Theta}).$ 
 Then, the following statements hold true.
\begin{itemize}
 \item[(a)] The stratifying system $(\Theta,\leq)$ is compatible with $I,$ and $(G\Theta,G\underline{Q},\leq)$ is the Ext-projective 
stratifying system associated to $(G\Theta,\leq).$
 \item[(b)] $G|_{\F(\Theta)}:\F(\Theta)\to \F(G\Theta)$ is an exact equivalence of categories and 
$F|_{\F(G\Theta)}:\F(G\Theta)\to \F(\Theta)$ is its quasi-inverse.
\end{itemize}
\end{teo}
\begin{dem} By \ref{GfielPleno} (a) we have that $\Hom_\Lambda(M,I\otimes N)=0$ for any $M,N\in\F(\Theta).$ Thus, by \ref{FielPlenoG}, we conclude 
that $G=G|_{\F(\Theta)}:\F(\Theta)\to \F(G\Theta)$ is a fully faithful functor. Furthermore, from \ref{GfielPleno} (b), it follows that 
$\Ext_\Lambda^1(\mathbf{\Theta},I\otimes\mathbf{\Theta})=0.$ In particular, the stratifying system
$(\Theta,\leq)$ is compatible with $I.$ Moreover, from \ref{G-epss}, we conclude that $(G\Theta,G\underline{Q},\leq)$ is the Ext-projective 
stratifying system associated to $(G\Theta,\leq);$ and hence (a) follows. 
\

In order to prove (b), it is enough to see that $G|_{\F(\Theta)}:\F(\Theta)\to \F(G\Theta)$ is dense. Indeed, if the restriction $G|_{\F(\Theta)}$ 
is a dense functor, we would have that $\Ima\,(G|_{\F(\Theta)})=\F(G\Theta);$ and so from \ref{Equivalencias} we conclude (b).
\

Finally, we prove that the functor $G=G|_{\F(\Theta)}:\F(\Theta)\to \F(G\Theta)$ is dense. Indeed, let 
$M\in\F(G\Theta).$ We proceed by induction on the $G\Theta$-length  
$\ell_{G\Theta}(M).$ If $\ell_{G\Theta}(M)=1$ then $M\simeq G\Theta(i)$ for some $i.$\\
Let $\ell_{G\Theta}(M)>1.$ Then, by \cite[Lemma 2.8]{MMS2}, there is an exact sequence 
$0\to G\Theta(i)\to M\to M/G\Theta(i)\to 0$ in $\modu\,(\Gamma),$ where
$\ell_{G\Theta}(M/G\Theta(i))=\ell_{G\Theta}(M)-1$ for some $i.$ So, by
induction, there exists $Z\in \F(\Theta)$ such that
$G(Z)\simeq M/G\Theta(i)$. Moreover, by \cite[Proposition 2.10]{MMS2}, there
is an exact sequence $\eta_{Z}\;:\;0\to Z'\stackrel{u}{\to} Q_{0}(Z)\stackrel{\varepsilon_{Z}}
{\to} Z\to 0$ in $\F(\Theta),$ with $Q_{0}(Z)\in\add\,(Q).$ Thus, we get the following 
exact and commutative diagram in $\modu\,(\Gamma)$
$$\xymatrix{ & & 0\ar[d] & 0\ar[d]\\
& & G(Z')\ar@{=}[r]\ar[d]^{\mu} & G(Z')\ar[d]^{G(u)}\\
\eta\;:\; 0\ar[r] & G\Theta(i)\ar[r]^{i_{1}}\ar@{=}[d] &
C\ar[r]^{p_2}\ar[d]^{\lambda} &
G(Q_{0}(Z))\ar[r]\ar[d]^{G(\varepsilon_{Z})} & 0\\
0\ar[r] & G\Theta(i)\ar[r] & M\ar[r]\ar[d] &
G(Z)\ar[r]\ar[d] & 0\\
& & 0 & 0.}$$ 
Since $G(Q_{0}(Z))$ is Ext-projective in $\F(G\Theta),$ the exact sequence $\eta$ splits and hence $C=\,
G\Theta(i)\bigoplus G(Q_{0}(Z))\simeq
G(\Theta(i)\bigoplus G(Z))$,
$i_{1}=\left(\begin{array}{c}
1\\
0
\end{array}\right)$ and $p_{2}=(0,1)$. That is $\mu=\left(\begin{array}{c}
\varphi\\
G(u) \end{array}\right)$ with $\varphi:G(Z') \to G(\Theta(i))$. Using that the restriction 
$G|_{\F(\Theta)}$ is full, there exists  $h:Z' \to \Theta(i)$ such that 
$G(h)=\varphi$ and hence $\mu=G(\psi),$ where $\psi:=\left(\begin{array}{c}
h\\
u
\end{array}\right)$. Observe that the morphism $\psi$ is a monomorphism since $u$ is so. Then, by completing $\psi$ to an 
exact sequence,  we get the following commutative diagram 
$$\xymatrix{ & & 0\ar[d] & 0\ar[d]\\
& & \Theta(i)\ar@{=}[r]\ar[d] &  \Theta(i)\ar[d]\\
0\ar[r] & Z'\ar[r]^(.4){\psi}\ar@{=}[d] & \Theta(i)\bigoplus
Q_{0}(Z)\ar[r]\ar[d]^{\pi_2} & X\ar[r]\ar[d]^{\alpha}
& 0\\
0\ar[r] & Z'\ar[r]^{u} & Q_{0}(Z)\ar[r]^{\varepsilon_{Z}}\ar[d] &
Z\ar[r]\ar[d] & 0\\
& & 0 & 0,}$$
where the rows and columns are exact sequences and $\pi_2$ is the canonical projection. Observe that 
$X\in \F(\Theta)$ since $\F(\Theta)$ is closed under
extensions. Thus, we get the exact sequence
$$\xymatrix{0\ar[r] & G(Z')\ar[r]^(.4){G(\psi)} & G(\Theta(i)\bigoplus
Q_{0}(Z))\ar[r] & G(X)\ar[r] & 0.}$$ But $G(\psi)=\mu$ and hence 
$G(X)\simeq \Coker\,(\mu)=M;$ proving that $G:\F(\Theta)\to \F(G\Theta)$ is dense.
\end{dem}
\vspace{.2cm}

Let us consider the following examples.

\begin{ex} Let $\Gamma$ be the quotient path $k$-algebra given by the quiver $$
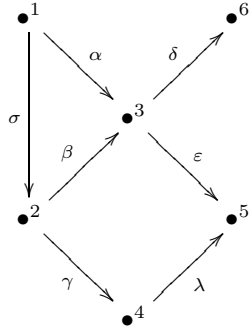
\xymatrix{\bullet^1\ar[dr]^{\alpha}\ar[dd]_{\sigma} & & \bullet^6\\
& \bullet^3\ar[ur]^{\delta}\ar[dr]^{\varepsilon} & \\
\bullet^2\ar[ur]^{\beta}\ar[dr]_{\gamma} & & \bullet^5\\
& \bullet^4\ar[ur]_{\lambda} & }$$
modulo the relations $\delta\beta=0$ and $\varepsilon\beta=\lambda\gamma.$ Consider the 
ideal $I=\la\overline{\beta},\overline{\gamma}\ra\unlhd \Gamma.$ So, $\Gamma$ is an 
split-by-nilpotent extension of $\Lambda:=\Gamma/I$ by $I.$ Furthermore, the algebra $\Lambda$ is the 
path $k$-algebra given by the quiver
$$
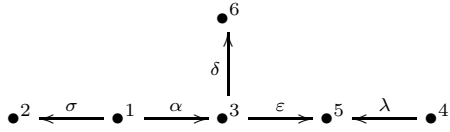
\xymatrix{&& \bullet^6 &&\\
\bullet^2& \bullet^1\ar[l]_{\sigma}\ar[r]^{\alpha} & \bullet^3\ar[u]^{\delta}\ar[r]^{\varepsilon} & \bullet^5 & \bullet^4\ar[l]_{\lambda}}$$
Observe that $I_\Lambda$ is projective, since $I_\Lambda\simeq {}_{\Lambda^{op}}P(2)^3.$ We consider the natural order $1\leq 2\leq 3\leq 4$ and the stratifying system $(\Theta,\leq)$ 
of size $4$ in $\modu\,(\Lambda),$ where 
$\Theta(1)=\begin{matrix}S(1)\\S(3)\end{matrix}={}_\Lambda I(3),$ $\Theta(2)=\begin{matrix}
S(3)\;S(4)\\  S(5) \end{matrix},$ $\Theta(3)=S(4),$ $\Theta(4)=S(2)={}_\Lambda P(2).$ An explicit 
calculation gives us that $I\otimes\Theta(i)=0$ for $i=1,2,3$ and $I\otimes\Theta(4)\simeq \Theta(2).$ So $I\otimes\mathbf{\Theta}\in\F(\Theta)$ and by \ref{FilItT} it follows that 
$\Ext_\Lambda^1(Q, I\otimes\mathbf{\Theta})=0.$ Moreover, it can be seen that the stratifying system $(\Theta,\leq)$ is compatible with the ideal $I.$ Thus, the needed conditions in \ref{GbarraInducido} hold. Finally, observe that $\Hom_\Lambda(\Theta,I\otimes\mathbf{\Theta})\neq 0$ since $\Hom_\Lambda(\Theta(2),I\otimes
\Theta(4))\simeq\End_\Lambda(\Theta(2))\neq 0.$
\end{ex}

\begin{ex} Let $\Gamma$ be the quotient path $k$-algebra given by the quiver $$
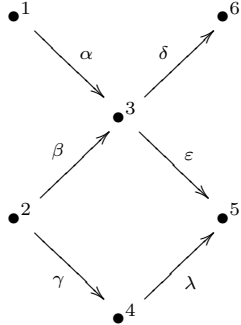
\xymatrix{\bullet^1\ar[dr]^{\alpha} & & \bullet^6\\
& \bullet^3\ar[ur]^{\delta}\ar[dr]^{\varepsilon} & \\
\bullet^2\ar[ur]^{\beta}\ar[dr]_{\gamma} & & \bullet^5\\
& \bullet^4\ar[ur]_{\lambda} & }$$
modulo the relations $\delta\beta=0$ and $\varepsilon\beta=\lambda\gamma.$ Consider the 
ideal $I=\la\overline{\beta},\overline{\gamma}\ra\unlhd \Gamma.$ So, $\Gamma$ is an 
split-by-nilpotent extension of $\Lambda:=\Gamma/I$ by $I.$ Furthermore, the algebra $\Lambda$ is the 
path $k$-algebra given by the quiver
$$
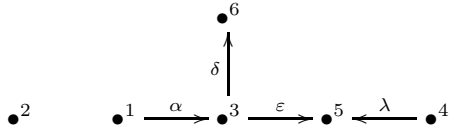
\xymatrix{&& \bullet^6 &&\\
\bullet^2&  \bullet^1\ar[r]^{\alpha} & \bullet^3\ar[u]^{\delta}\ar[r]^{\varepsilon} & \bullet^5 & \bullet^4\ar[l]_{\lambda}}$$
Observe that $I_\Lambda$ is projective, since $I_\Lambda\simeq {}_{\Lambda^{op}}P(2)^3.$ We consider the natural order $1\leq 2\leq 3$ and the stratifying system $(\Theta,\leq)$ 
of size $3$ in $\modu\,(\Lambda),$ where 
$\Theta(1)=\begin{matrix}\quad\quad\; S(1)\\S(4)\;S(3)\\S(5)\end{matrix}={}_\Lambda I(5),$ $\Theta(2)=S(2)={}_\Lambda P(2)$ and $\Theta(3)=S(6)={}_\Lambda P(6).$ An explicit 
calculation gives us that $I\otimes\Theta(i)=0$ for $i=1,3$ and $N:=I\otimes\Theta(2)=\begin{matrix}
S(3)\;S(4)\\  S(5) \end{matrix}.$ It can be seen that  
$\Hom_\Lambda(\mathbf{\Theta}, I\otimes\mathbf{\Theta})=0=\Ext_\Lambda^1(\mathbf{\Theta}, I\otimes\mathbf{\Theta}).$ Thus, by \ref{GfielPleno}, 
the needed conditions in \ref{Gequiv} hold. Finally, observe that $I\otimes\mathbf{\Theta}\not\in \F(\Theta)$ since $I\otimes
\Theta(2)=N\not\in\F(\Theta).$
\end{ex}

We finish this section by taking into consideration the canonical stratifying system $({}_\Lambda\Delta,\leq ).$ 
Let the standard 
$\Lambda$-modules ${}_\Lambda\Delta$ be computed using the pair $({}_\Lambda P,\leq),$ where 
${}_\Lambda P=\{{}_\Lambda P(i)\}_{i=1}^n$ is a representative set of the indecomposable projective $\Lambda$-modules, $n:=rk\,K_0(\Lambda)$ 
and $\leq$ is a linear order on the set $[1,n].$ By \ref{EquiProj}, we have that ${}_\Gamma P:=G({}_\Lambda P)$ 
is a representative set of the indecomposable projective $\Gamma$-modules. So, we compute the standard $\Gamma$-modules ${}_\Gamma\Delta$ 
by using the pair $({}_\Gamma P,\leq).$

\begin{teo}\label{ApliGss-alg} Let $I_{\Lambda}$ be projective, and let $({}_\Lambda\Delta,\leq)$ be compatible with the ideal $I.$ Then, the 
following statements hold true.  
\begin{itemize}
 \item[(a)]  $G({}_\Lambda\Delta(i))\simeq {}_\Gamma\Delta(i)$ for any $i\in[1,n].$ 
 \item[(b)]  $G|_{\F({}_\Lambda\Delta)}:\F({}_\Lambda\Delta)\to \F({}_{\Gamma}\Delta)$ is an exact and faithful functor which preserves and reflects indecomposables.
 \item[(c)] $F|_{\F({}_\Gamma\Delta)}:\F({}_\Gamma\Delta)\to \F({}_{\Lambda}\Delta)$ is an exact, faithful and dense functor which reflects indecomposables.  
\end{itemize}
\end{teo}

\begin{dem} By \ref{FiltrFG}, we have the    
stratifying system $(G({}_\Lambda\Delta),\leq)$ in $\modu\,(\Gamma)$ and also the exact functors  $\F({}_\Lambda\Delta)\stackrel{G}{\to} \F(G({}_{\Lambda}\Delta))\stackrel{F}{\to}\F({}_\Lambda\Delta)$ satisfying the desired properties as in (b) and (c). It remains to show that $G({}_\Lambda\Delta(i))\simeq {}_\Gamma\Delta(i)$ for any $i\in[1,n].$ Indeed, for each 
$i\in[1,n],$ consider the exact sequence 
$$\eta_i\;:\;0\longrightarrow K(i)\longrightarrow {}_\Lambda P(i)\stackrel{p_i}{\longrightarrow}{}_\Lambda\Delta(i)\longrightarrow 0,$$
\noindent where $K(i):=\Tr_{\oplus_{j>i}\,{}_\Lambda P(j)}\,({}_\Lambda P(i)).$ Thus, we get the exact sequence 
$$G(\eta_i)\;:\;0\longrightarrow G(K(i))\longrightarrow {}_\Gamma P(i)\stackrel{G(p_i)}{\longrightarrow}G({}_\Lambda\Delta(i))\longrightarrow 0.$$
\noindent Let $Z(i):=\Tr_{\oplus_{j>i}\,{}_\Gamma P(j)}\,({}_\Gamma P(i)),$  and let $\alpha:{}_\Gamma P(j)\to {}_\Gamma P(i),$ with $j>i$ be any morphism. Then, from \ref{EquiProj}, there exists $\alpha':{}_\Lambda P(j)\to {}_\Lambda P(i)$ such that $G(
\alpha')=\alpha;$ and hence $G(p_i)\alpha=G(p_i\alpha')=0.$ Therefore $\Ima\,(\alpha)\subseteq 
G(K(i));$ proving that $Z(i)\subseteq G(K(i)).$ On the other hand, since $K(i):=\Tr_{\oplus_{j>i}\,{}_\Lambda P(j)}\,({}_\Lambda P(i))$ and $G$ is an exact functor, we 
get an epimorphism  $\oplus_{j>i}\,{}_\Gamma P(j)^m\to G(K(i)),$ getting us that $G(K(i))
\subseteq Z(i).$  
\end{dem}

\begin{cor}\label{ApliGss-algSS} Let $I_{\Lambda}$ be projective, and let $({}_\Lambda\Delta,\leq)$ be compatible with the ideal $I.$ Then the 
following statements hold true.  
\begin{itemize}
\item[(a)]  $\Lambda$ is an standardly stratified (respectively, quasi-hereditary) algebra if and only if $\Gamma$ is so.
 \item[(b)] Let $\Lambda$ be an standardly stratified algebra such that $\Hom_\Lambda(\mathbf{{}_\Lambda\Delta},I\otimes\mathbf{{}_\Lambda\Delta})=0.$ Then, the functor $G|_{\F({}_\Lambda\Delta)}:\F({}_\Lambda\Delta)\to \F({}_{\Gamma}\Delta)$  is an exact equivalence with quasi-inverse $F|_{\F({}_\Gamma\Delta)}:\F({}_\Gamma\Delta)\to \F({}_{\Lambda}\Delta).$
\end{itemize}
\end{cor}
\begin{dem}
(a) Let $\Lambda$ be  a standardly stratified algebra. For each 
$i\in[1,n],$ consider the exact sequence 
$$\eta_i\;:\;0\longrightarrow K(i)\longrightarrow {}_\Lambda P(i)\stackrel{p_i}{\longrightarrow}{}_\Lambda\Delta(i)\longrightarrow 0,$$
\noindent where $K(i):=\Tr_{\oplus_{j>i}\,{}_\Lambda P(j)}\,({}_\Lambda P(i)).$ Since $\Lambda$ is an ss-algebra, we have that $\eta_i$ lies in 
$\F({}_\Lambda\Delta).$ Thus, by \ref{ApliGss-alg} (a), we conclude that the exact sequence $G(\eta_i)$ lies in $\F({}_\Gamma\Delta);$ 
proving that $\proj\,(\Gamma)\subseteq \F({}_\Gamma\Delta).$ That is, the algebra $\Gamma$ is 
standardly stratified.
\

Assume now that $\Gamma$ is standardly stratified. Then, by \cite{AHLU}, it follows that 
$\F({}_\Gamma\Delta)$ is a resolving category. In particular, the exact sequence $G(\eta_i)$ 
lies in $\F({}_\Gamma\Delta).$ Since $F|_{\F({}_\Gamma\Delta)}:\F({}_\Gamma\Delta)
\to \F({}_{\Lambda}\Delta)$ is exact, by applying $F$ to $G(\eta_i),$ we get that the exact 
sequence $\eta_i$ lies in $\F({}_\Lambda\Delta)$ (see \ref{basicos2} (a));  proving that 
$\proj\,(\Lambda)\subseteq \F({}_\Lambda\Delta).$ That is, the algebra 
$\Lambda$ is standardly stratified.    
\

The proof that $\Lambda$ is quasi-hereditary if and only if $\Gamma$ is so, can be done as we 
did in the proof of \ref{GbarraInducido} (d). To do so, just replace there $\Theta$ by 
${}_\Lambda\Delta$ and use \ref{ApliGss-alg} (a).
\

(b) Let $\Lambda$ be a standardly stratified algebra. Then, the triple $({}_\Lambda\Delta,{}_\Lambda P,\leq)$ is the Ext-projective stratifying system associated to $({}_\Lambda\Delta,\leq).$ So the result follows from \ref{Gequiv}.
\end{dem}

\section{Restricting stratifying systems with the functor $F$}

In this section, let $\Gamma$ be a split-by-nilpotent extension of $\Lambda$ by and ideal 
$I\unlhd\Gamma.$ We consider the change of rings functor $F:={}_{\Lambda}\Lambda_{\Gamma} \otimes_{\Gamma} -:
\modu\,(\Gamma)\to\modu\,(\Lambda).$ We ask under which conditions a stratifying system in $\modu\,(\Gamma)$ can be restricted, through the functor $F,$ to a stratifying system in $\modu\,(\Lambda).$
\

Observe that, by \ref{basicos2} (d), we know that $F$ reflects indecomposables. However, the functor $F,$ in general, could not preserve indecomposables. A sufficient condition, for the functor $F,$ to preserve indecomposability is given in the following proposition. In order to do that, we will need the next lemma.

\begin{lem}\label{basicoF} Let $\Lambda_\Gamma$ be projective. Then, for any $M,N \in \modu\,(\Gamma)$, we have the following long exact sequence of $R$-modules
\[ 0\rightarrow \Hom_{\Gamma}(M, I\otimes_{\Gamma} N) \rightarrow \Hom_{\Gamma}(M,N)\rightarrow \Hom_{\Lambda}(F(M),F(N)) \rightarrow \]
\[\Ext_{\Gamma}^{1}  (M, I\otimes_{\Gamma} N) \rightarrow \Ext_{\Gamma}^{1}  (M, N) \rightarrow \Ext_{\Lambda}^{1}  (F(M), F(N)) \rightarrow \Ext_{\Gamma}^{2}  (M, I\otimes_{\Gamma} N). \]
\end{lem}
\begin{dem} Let $M,N \in \modu\,(\Gamma).$ We assert that 
\[ \Ext_{\Gamma}^{i}  (M, \Lambda \otimes_{\Gamma} N) \simeq \Ext_{\Lambda}^{i}  (F(M),F(N)) \ \text{ for any} \ i.\]
Indeed, we get that $\Ext_{\Lambda}^{i}  (F(M), F(N))\simeq \Ext_{\Gamma}^{i}  (M, \Hom_{\Lambda} ({}_{\Lambda}\Lambda_{\Gamma}, F(N)))$ since $\Lambda_\Gamma$ is projective (see \cite[Exercise 9.21]{Ro}). On the other hand,  from \ref{CambioAni2},
we have that $\Hom_{\Lambda} ({}_{\Lambda}\Lambda_{\Gamma},{}_{\Lambda}\Lambda_{\Gamma} \otimes_{\Gamma} N) \simeq \Hom_{\Lambda} ({}_{\Lambda}\Lambda_{\Gamma}, {}_{\Lambda}\Lambda_{\Gamma}) \otimes_{\Gamma} N\simeq{}_\Gamma\Lambda_\Gamma\otimes_{\Gamma} N$ as $\Gamma$ modules; and so the assertion follows.
\

Applying the functor $- \otimes_{\Gamma} N$ to the exact sequence of  $\Gamma-\Gamma$ bimodules $0 \rightarrow I \rightarrow \Gamma \xrightarrow{\pi} \Lambda \rightarrow 0$, and using the fact that $\Tor_{1}^{\Gamma} (\Lambda_{\Gamma},\ _{\Gamma}N)=0$, we get the exact sequence of  $\Gamma$-modules
$\eta\;:\; 0 \rightarrow I\otimes_{\Gamma}N \rightarrow N \rightarrow \Lambda \otimes_{\Gamma} N \rightarrow 0.$ Moreover, by applying the functor $\Hom_{\Gamma} (M,-)$ to $\eta$, and the above 
assertion, we get the desired exact sequence. 
\end{dem}

\begin{pro}\label{F-indesc} Let $\Lambda_\Gamma$ be projective and $M\in\modu\,(\Gamma)$ be indecomposable. If 
$\Hom_\Gamma(M,I\otimes_\Gamma M)=0=\Ext_\Gamma^1(M,I\otimes_\Gamma M)$ then $F(M)$ is 
indecomposable.
\end{pro}
\begin{dem} By the assumed conditions and \ref{basicoF}, we get an isomorphism 
 $\xi:\End_\Gamma(M)\to \End_\Lambda(F(M))$ of $R$-modules. An explicit computation of the 
 map $\xi$ gives us that $\xi(f)(\lambda\otimes m)=\lambda\otimes f(m)$ for any $\lambda\in\Lambda,$ $m\in M$ and $f\in \End_\Gamma(M).$ Thus, the map $\xi$ is also a ring homomorphism and hence $\End_\Gamma(M)\simeq \End_\Lambda(F(M))$ as $R$-algebras; proving that 
 $F(M)$ is indecomposable. 
\end{dem}

\begin{defi} A stratifying system $(\Psi,\leq),$ of size $t$ in $\modu\,(\Gamma),$ is {\bf 
compatible} with the ideal $I\unlhd \Gamma$ if the following conditions hold.
\begin{itemize}
\item[(C1)] $\Ext_\Gamma^1(\Psi(j),I\otimes_{\Gamma}\Psi(i))=0$ for $j>i.$ 
\vspace{.2cm}
\item[(C2)] $\Ext_\Gamma^2(\Psi(j),I\otimes_{\Gamma}\Psi(i))=0$ for $j\geq i.$
\end{itemize}
\end{defi}

\begin{teo}\label{TeoRestF} Let $\Lambda_\Gamma$ be projective, and let $(\Psi,\leq)$ be a 
stratifying system of size $t$ in $\modu\,(\Gamma),$ which is compatible with the ideal $I\unlhd \Gamma.$ Then, for each $i\in[1,t]$ and any choice of an indecomposable direct summand $\Theta(i)$ of $F(\Psi(i)),$ 
the pair $(\Theta,\leq)$ is a stratifying system of size $t$ in $\modu\,(\Lambda).$  
\end{teo}
\begin{dem} Let $i\in[1,t].$ Since $\Psi(i)\neq 0,$ it follows from \ref{basicos2} (b) that 
$F(\Psi(i))\neq 0.$ Thus $F(\Psi(i))$ has at least one indecomposable direct summand. Let $\Theta(i)$ be a choice of an indecomposable direct summand of $F(\Psi(i)).$ Since $(\Psi,\leq)$ 
is compatible with the ideal $I,$ we get from \ref{basicoF} that the following conditions hold: (a) $\Hom_\Lambda(F(\Psi(j)),F(\Psi(i)))=0$ for $j>i,$ and (b) $\Ext_\Lambda^1(F(\Psi(j)),F(\Psi(i)))=0$ for $j\geq i.$ Therefore, the same conditions, as in (a) and (b), hold for $\Theta:=\{\Theta(i)\}_{i=1}^t$ since each 
$\Theta(i)$ is an indecomposable direct summand of $F(\Psi(i)).$ Hence the result follows.
\end{dem}

\begin{ex} Let $\Gamma$ be the split-by-nilpotent extension  of $\Lambda$ by $I,$ which is considered in \ref{Ej2}, and take the natural order $1\leq 2\leq 3.$ Consider the pair $(\Psi,\leq),$ where $\Psi(1):={}_\Gamma S(2),$ $\Psi(2):={}_\Gamma P(1)$ and $\Psi(3):={}_\Gamma S(3).$
Observe that $(\Psi,\leq)$ is a stratifying system in $\modu\,(\Gamma).$ Furthermore, since 
$I\otimes_\Gamma\Psi(1)=0=I\otimes_\Gamma\Psi(2),$ $I\otimes_\Gamma\Psi(3)\simeq {}_\Gamma S(1)$ and $\id\,({}_\Gamma S(1))\leq 1,$ it can be seen that the pair $(\Psi,\leq)$ is compatible with 
the ideal $I;$ and furthermore $\Lambda_\Gamma$ is projective since $\Lambda_\Gamma\simeq 
{}_{\Gamma^{op}}P(3)\oplus {}_{\Gamma^{op}}P(2).$  Thus by \ref{TeoRestF} it follows that, for any choice of an indecomposable direct summand $\Theta(i)$ of $F(\Psi(i)),$ 
the pair $(\Theta,\leq)$ is a stratifying system of size $3$ in $\modu\,(\Lambda).$ Since 
$F(\Psi(1))\simeq {}_\Lambda P(2),$ $F(\Psi(2))\simeq {}_\Lambda P(1)$  and $F(\Psi(3))\simeq {}_\Lambda P(3),$ we see that the restriction of $(\Psi,\leq)$ to $\modu\,(\Lambda),$ through the functor $F,$ gives only one stratifying system.
\end{ex}

\footnotesize

\vskip3mm \noindent Marcelo Lanzilotta M:\\
Centro de Matem\' atica (CMAT), \\ 
Instituto de Matem\'atica y Estad\'\i stica Rafael Laguardia (IMERL),\\
Universidad de la Rep\'ublica.\\ Igu\'a 4225, C.P. 11400, Montevideo, URUGUAY.

{\tt marclan@cmat.edu.uy, marclan@fing.edu.uy}

\vskip3mm \noindent Octavio Mendoza:\\
Instituto de Matem\'aticas,\\
Universidad Nacional Aut\'onoma de M\'exico,\\
Circuito Exterior, Ciudad Universitaria,\\
M\'exico D.F. 04510, M\'EXICO.

{\tt omendoza@matem.unam.mx}

\vskip3mm \noindent Corina S\'aenz:\\
Departamento de Matem\'aticas, Facultad de Ciencias,\\
Universidad Nacional Aut\'onoma de M\'exico,\\
Circuito Exterior, Ciudad Universitaria,\\
M\'exico D.F. 04510, M\'EXICO.

{\tt ecsv@lya.fciencias.unam.mx}

\end{document}